\documentclass[12pt]{article}
\usepackage{auto-pst-pdf} % use with pstricks Use Tex & Dvi
\usepackage{pstricks,pst-node}
\usepackage{amsmath,amsthm}
\usepackage{amssymb}
\usepackage{color}
\usepackage{xspace}
\usepackage[colorlinks=true,
linkcolor=green,
filecolor=brown,
citecolor=green]{hyperref}

\def\s{Schr\"{o}der\xspace}
\def\gf{generating function\xspace}

\begin{document}
\newtheorem*{prop*}{Proposition}
\newtheorem{theorem}{Theorem}
\newtheorem{defn}[theorem]{Definition}
\newtheorem{lemma}[theorem]{Lemma}
\newtheorem{prop}[theorem]{Proposition}
\newtheorem{cor}[theorem]{Corollary}
%\vspace*{-5mm}
\begin{center}
{\Large
A note on a bijection for Schr\"{o}der permutations                           \\ 
}
\vspace{10mm}
DAVID CALLAN  \\
Department of Statistics  \\
\vspace*{-1mm}
University of Wisconsin-Madison  \\
\vspace*{-1mm}
Madison, WI \ 53706-1532  \\
{\bf callan@stat.wisc.edu}  \\
\vspace{5mm}

\end{center}

\begin{abstract}
There is a bijection from Schr\"{o}der paths to $\{4132,\,4231\}$-avoiding permutations due to Bandlow, Egge, and Killpatrick that sends ``area'' to ``inversion number''. Here we give a concise description of this bijection.
\end{abstract}

\section{Introduction}

\vspace*{-3mm}

Permutations of length $n$ avoiding the two 4-letter patterns 4132 and 4231 are 
known \cite{kremer2000} to be counted by the large 
Schr\"{o}der number $r_{n-1}$  with \gf
$\sum_{n\ge 0}r_n x^n =1 +2x+6x^2+22x^3+\dots = \frac{1 - x - \sqrt{1 - 6x + x^2}}{2x}$ 
(\htmladdnormallink{A006318}{http://oeis.org/A006318}).
It is well known that $r_n$ is the number of Schr\"{o}der paths of size $n$ where a \emph{Schr\"{o}der path} is a lattice path of north steps $N=(0,1)$, 
diagonal steps $D=(1,1)$ and east steps $E=(1,0)$ 
that starts at the origin, never drops below the diagonal $y=x$, and terminates on the diagonal. 
Its \emph{size} is $\#\,N$ steps + $\#\,D$ steps, and a  
Schr\"{o}der $n$-path is one of size $n$. Thus a Schr\"{o}der $n$-path ends at $(n,n)$. 

Bandlow, Egge, and Killpatrick \cite{bandlow2002} construct a bijection from Schr\"{o}der $n$-paths to \{4132, 4231\}-avoiding permutations of length $n+1$. They show that their bijection relates the area below the path to the number of inversions in the permutation. Here we give a concise, and possibly more transparent, formulation of their bijection.

It is convenient, and natural for present purposes, to take our permutations on $[0,n]$ rather than the standard $[1,n]$ and to define
the \emph{truncated coinversion table} of a permutation $p$ on $[0,n]$ to be the list $c=(c_i)_{i=1}^{n}$ 
where $c_i$ is the number of coinversions topped by $i$, that is, $c_i=\#\,\{j:0\le j<i 
$ and $j$ precedes $i$ in the permutation\}. 
The table is ``truncated'' because we omit $c_0$ which would necessarily be 0.
We say $i$ is a fixed point of $c$ if $c_i=i$  and denote the list of fixed points of $c$ by FP($c$).
Thus, for $p=350214,$ we have $c=11041$ and FP($c)=(1,4)$. Clearly, the fixed points of $c$ are right-to-left minima of $p$ and conversely, every nonzero right-to-left minimum of $p$ is a fixed point of $c$.

\section{The Bandlow-Egge-Killpatrick bijection}

\vspace*{-3mm}

In these terms, the bijection from \{4132, 4231\}-avoiders $p$ on $[0,n]$ to \s $n$-paths 
has a short description: find the coinversion table $c$ of $p$, then form the unique 
Schr\"{o}der $n$-path whose $D$ steps end at $x$ coordinates in FP($c$) and whose $N$ 
steps, taken in order, have $x$ coordinates given by the list $c\, \setminus \,$FP($c$).
For example, with $n=10$ and $p=2\ 1\ 0\ 7\ 9\ 6\ 10\ 5\ 3\ 4\ 8$, we have 
$c=0\ 0\ 3\ 4\ 3\ 3\ 3\ 8\ 4\ 6$, 
FP($c)=3\ 4\ 8$, and $c\, \setminus \,$FP($c)=0\ 0\ 3\ 3\ 3\ 4\ 6$.
The corresponding Schr\"{o}der path is shown in Figure 1.

\begin{center}
\begin{pspicture}(0,-1.5)(6,6.5)%\showgrid
\psset{unit=.6cm}
\newrgbcolor{mycolor}{1 1 0.2}
\pspolygon[fillstyle=solid,fillcolor=mycolor,linecolor=lightgray](0,0)(0,1)(0,2)(1,2)(2,2)(3,2)(3,3)(3,4)(3,5)(3,6)(4,6)(4,7)(4,8)(5,8)(6,8)(6,9)(7,9)(8,9)(8,10)(9,10)(10,10)(0,10)(0,0)
\psset{linecolor=lightgray}
\psdots(0,0)(0,1)(0,2)(0,3)(0,4)(0,5)(0,6)(0,7)(0,8)(0,9)(0,10)
(1,1)(1,2)(1,3)(1,4)(1,5)(1,6)(1,7)(1,8)(1,9)(1,10)
(2,2)(2,3)(2,4)(2,5)(2,6)(2,7)(2,8)(2,9)(2,10)
(3,2)(3,3)(3,4)(3,5)(3,6)(3,7)(3,8)(3,9)(3,10)
(4,4)(4,5)(4,6)(4,7)(4,8)(4,9)(4,10)
(5,5)(5,6)(5,7)(5,8)(5,9)(5,10)
(6,6)(6,7)(6,8)(6,9)(6,10)
(7,7)(7,8)(7,9)(7,10)(8,8)(8,9)(8,10)(9,9)(9,10)(10,10)
\psline(0,0)(0,10)(10,10)
\psline[linestyle=dashed](0,0)(10,10)
\psset{linecolor=black}
\psdots(0,0)(0,1)(0,2)(1,2)(2,2)(3,3)(3,4)(3,5)(3,6)(4,7)(4,8)(5,8)(6,8)(6,9)(7,9)(8,10)(9,10)(10,10)
\psline[linewidth=1.2pt,linecolor=black](0,0)(0,1)(0,2)(1,2)(2,2)(3,3)(3,4)(3,5)(3,6)(4,7)(4,8)(5,8)(6,8)(6,9)(7,9)(8,10)(9,10)(10,10)

\rput(0.2,-.3){\textrm{{\scriptsize $0$}}}
\rput(1.2,.7){\textrm{{\scriptsize $1$}}}
\rput(2.2,1.7){\textrm{{\scriptsize $2$}}}
\rput(3.2,2.7){\textrm{{\scriptsize $3$}}}
\rput(4.2,3.7){\textrm{{\scriptsize $4$}}}
\rput(5.2,4.7){\textrm{{\scriptsize $5$}}}
\rput(6.2,5.7){\textrm{{\scriptsize $6$}}}
\rput(7.2,6.7){\textrm{{\scriptsize $7$}}}
\rput(8.2,7.7){\textrm{{\scriptsize $8$}}}
\rput(9.2,8.7){\textrm{{\scriptsize $9$}}}
\rput(10.2,9.7){\textrm{{\scriptsize $10$}}}

\rput(5,-1.2){\textrm{A Schr\"{o}der 10\hspace*{1pt}-path }}
\rput(5,-2.5){\textrm{Figure 1}}

\end{pspicture}
\end{center}

\noindent The map is both defined and reversible due to the following characterization of $\{4132, 4231\}$-avoiders whose straightforward proof is left to the interested reader.
\begin{prop*}
Suppose $p$ is a permutation on $[0,n]$ with truncated coinversion table $c$. Then $p$ avoids $\{4132, 4231\}$ if and only if the list $c\, \setminus \,$\emph{FP}\hspace*{.5pt}$(c)$ is weakly increasing.  \qed
\end{prop*}

The connection between number of inversions in the permutation (= $\binom{n+1}{2}-\#$ coinversions) and area of the path is clear: the number of unit squares in the first quadrant lying wholly or partially to the left of the path (in yellow in Figure 1) is the sum of the $x$ coordinates of all the north and diagonal steps which, by the construction, is $\sum_{i=1}^{n}c_i$, the total number of coinversions in the permutation.

\end{document}